\documentclass[11pt]{amsart}
\usepackage{graphics,verbatim,color}

\usepackage[all]{xy}

\newtheorem{theorem}{Theorem}
\newtheorem*{theorem*}{Theorem}

\newtheorem*{proposition*}{Proposition}

\newtheorem*{definition*}{Definition}

\newtheorem*{conjecture*}{Conjecture}

\newtheorem*{problem*}{Problem}

\newtheorem*{example*}{Example}

\def\R{\mathbb{R}}

\def\C{\mathbb{C}}

\def\Q{\mathbb{Q}}

\def\g{\gamma}

\def\w{\omega}

\def\z{\zeta}

\def\smallskip{\par\vspace{1mm}}
\def\medskip{\par\vspace{2mm}}
\def\bigskip{\par\vspace{3mm}}

\def\fr#1#2{\frac{#1}{#2}}

\newcount\equationcount
\def\thenumber{0}
\def\eq#1{\global\advance\equationcount by 1
   \def\thenumber{\number\equationcount}
                        {$$#1\eqno(\thenumber)$$}}

\begin{document}

\title[Partially rigid motions]{Partially rigid motions in the planar three-body problem}
\author{Richard Moeckel}
\address{School of Mathematics\\ University of Minnesota\\ Minneapolis MN 55455}

\email{rick@math.umn.edu}



\begin{abstract}
A solution of the $n$-body problem in $\R^d$ is a {\em relative equilibrium} if all of the mutual distance between the bodies are constant.  In other words, the bodies undergo a rigid motion.  Here we investigate the possibility of partially rigid motions, where some but not all of the distances are constant.  For the planar three-body problem with equal masses, we show that partially rigid motions are impossible  -- if even one of the three mutual distances is constant, the motion must be a relative equilibrium. 
\end{abstract}

\date{May 9, 2025}
\maketitle

\section{Introduction}
Lagrange  discovered the {\em relative equilibrium} solutions of the three-body problem in $\R^3$, characterized by the  property that  all three mutual distances remain constant.   Up to rotation and scaling there are exactly four possible shapes for a relative equilibrium configuration -- the equilateral triangle configuration and three collinear configurations, one for each choice of which mass is between the other two.  

In this paper, we consider solutions which are only partially rigid.
\begin{definition*}
A solution of the $n$-body problem will be called {\em  partially rigid} if at least one, but not all, of the mutual distances remains constant.
\end{definition*}
In a previous paper we considered ``hinged'' solutions where all but one of the mutual distances.  It turns out that such solutions don't exist for $n=3,4$:
\begin{theorem*}[RM\cite{Moe}]\label{th_unhinged}
The three- and four-body problems in $\R^d$ are unhinged.  In other words, for $n=3$, the only solutions such that two of the mutual distances remain constant are the relative equilibrium solutions.  For $n=4$, the only solutions such that five of the six distance are constant are the relative equilibria.
\end{theorem*}
Note that for $n=3$, the theorem is false if two of the masses are zero.  If $m_1=m_2=0$ we can put them on circular orbits around $m_3$ to produce hinged solutions.  In fact this is the standard model of a clockwork solar system with two planets.  The theorem shows that such circular motions are impossible if the masses are all positive.

It seems plausible that fixing even one distances implies complete rigidity.
\begin{conjecture*}
There are no partially rigid solutions of the $n$-body problem in $\R^d$.  In other words, the only solutions such that one of the mutual distances remains constant are the relative equilibrium solutions.
\end{conjecture*}
The conjecture is false if even one of the masses is zero.  In fact, for the well-known circular restricted  three-body problem with $m_3=0$, all of the solutions have $r_{12}(t)=1$.   

The goal of this paper is to prove the full conjecture for the special case of the equal mass $3$-body problem in the plane, that is, the case $n=3$ and $d=2$.
\begin{theorem}\label{th_nopartialrigidity}
There are no partially rigid solutions of the planar three-body problem with equal masses.  In other words, the only solutions for which one or more of the three mutual distances remain constant are the relative equilibrium solutions of Lagrange.
\end{theorem}
In \cite{Moe} the elegant reduced differential equations for the $n$-body problem in $\R^d$ were used \cite{Alb, AC, Ch}.  Here we use Jacobi variables in the plane which are less elegant but are simpler and involve fewer variables.   The proof eventually reduces to showing that a system of five equations in two variable has only finitely many solutions.  In principle, this should be relatively straightforward but the equations are so complicated that an indirect argument using computer assisted algebraic computation with Newton polygons was needed.  For the three-body problem with unequal masses or in dimensions greater than two, even the reduction to two unknowns was not feasible.

The motivation for this problem came from considering a simple model of tidal friction.  Imagine that two of the three bodies are connected by a massless spring with friction.   When the length of the spring is changing, energy is dissipated and it is natural to conjecture that the system will eventually converge to some kind of equilibrium state where the length of the spring remains constant.   Is such a state necessarily a relative equilibrium or could there still be relative motions among the bodies~?  The conjecture about lack of partially rigid solutions is just the analogous problem without the spring.

\section{The planar three-body problem in Jacobi coordinates}
Consider the three-body problem with masses $m_i>0$ and positions $q_i\in\R^2$.  The problem is symmetric under translations and rotations.
We can eliminate translations by using Jacobi coordinates
$$z_1 = q_2-q_1\qquad z_2 = q_3-\nu_1q_1-\nu_2 q_2$$
where  $\nu_i = m_i/(m_1+m_2)$.  If $\z_i = \dot z_i$, the translation-reduced problem is a Lagrangian dynamical system with Lagrangian
$$L = \fr12(\mu_1 |\z_1|^2 + \mu_2 |\z_2|^2)+U(z_1,z_2)$$
where
$$U = \fr{m_1 m_2}{r_{12}} +  \fr{m_1 m_3}{r_{13}}+ \fr{m_2 m_3}{r_{23}},$$

$$r_{12} = |z_1|\qquad r_{13} = |z_2+\nu_2 z_1|\qquad r_{23} = |z_2-\nu_1 z_1|$$
and
$$\mu_1= \fr{m_1m_2}{m_1+m_2}\qquad \mu_2=\fr{(m_1+m_2)m_3}{m_1+m_2+m_3}.$$
The phase space is the tangent bundle of $\R^4\setminus \Delta$ where 
$$\Delta = \{(z_1,z_2)\in\R^4: r_{ij} = 0\,\text{for some }\,i\ne j\}.$$
In the special case of equal masses $m_1=m_2=m_3=1$ the mass parameters are
$$\mu_1=\fr12\quad \mu_2=\fr23\quad \nu_1=\nu_2=\fr12.$$

The problem is still symmetric under rotations of the $z_i, \z_i$ in the plane and the angular momentum is a constant of motion. Let
$$z_1=(x_1,y_1)\quad z_2=(x_2,y_2)\quad \z_1=(u_1,v_1)\quad \z_2=(u_2,v_2).$$
Then the angular momentum is
$$\Omega(z,\z) = \mu_1(x_1v_1-y_1u_1)+\mu_2(x_2v_2-y_2u_2).$$
The total energy
$$H(z,\z)=  \fr12(\mu_1 |\z_1|^2 + \mu_2 |\z_2|^2) - U(z_1,z_2)$$
is also constant.

\section{Solutions with a constant distance}
Suppose $\g(t)$ is a solution with $r_{12}(t)$ constant.  We will show that the other two distances $r_{13}(t), r_{23}(t)$ must also be constant which implies that the solution is a relative equilibrium.  Because of scaling symmetry, we may assume $r_{12}(t)=1$ for all $t$.  Clearly, all of the time derivatives of $f(t) = \fr12 r_{ij}(t)^2$ must vanish.   Each derivative imposes a constraint on the possible values of the eight coordinate variables $x_1,y_1,x_2,y_2,u_1,v_1,u_2,v_2$.   The strategy of the proof is to show that these equations, together with the energy and angular momentum constraints, restrict $r_{13}, r_{23}$ to a finite set.  Then by continuity, $r_{13}(t), r_{23}(t)$ must  be constant.

We have $f(t) = \fr12|z_1|^2 = \fr12(x_1^2+y_1^2)=\fr12$ and the first derivative is $f_1 = x_1u_1+y_1v_1$.  The second derivative is
$$f_2= u_1^2+v_1^2+x_1\dot u_1+y_1\dot v_1$$
where the velocity derivatives are computed using the Euler-Lagrange equations.  For the equal-mass case these are
$$\begin{aligned}
\dot u_1 &=  -\fr{2 x_1}{r_{12}^3} -\fr{x_2+2x_1}{r_{13}^3}-\fr{x_2-2x_1}{r_{23}^3}\\
\dot v_1 &=  -\fr{2 y_1}{r_{12}^3} -\fr{y_2+2y_1}{r_{13}^3}-\fr{y_2-2y_1}{r_{23}^3}\\
\dot u_2 &=  -\fr{3(x_2+2x_1)}{4r_{13}^3}-\fr{3(x_2-2x_1)}{4r_{23}^3}\\
\dot u_2 &=  -\fr{3(y_2+2y_1)}{4r_{13}^3}-\fr{3(y_2-2y_1)}{4r_{23}^3}.
\end{aligned}
$$
Similarly, using Mathematica, we can explicitly compute the higher time-derivatives as rational functions $f_k$, $k=3,4$, of the eight coordinate variables and the distance $r_{ij}$.

Relating these eleven variables, we have the equations $r_{12}=1$ and  $f_1=f_2=f_3=f_4=0$, together with formulas for the squared distances
$$r_{12}^2 = x_1^2+y_1^2  \quad r_{13}^2 =(x_2+\fr12 x_1)^2+ (y_2+\fr12 y_1)^2  \quad r_{23}^2 =(x_2-\fr12 x_1)^2+ (y_2-\fr12 y_1)^2 $$
and the energy and angular momentum equations
$$\Omega(z,\z)=\w \quad H(z,\z)=h$$
where $h,\w$ are constant parameters. 

We have ten equations in eleven variables and we can also use the rotational symmetry.   Assume, without loss of generality that at time $t=t_0$, $z_1$ lies along the positive $x$-axis.  Then at time $t=t_0$ we can eliminate three of the variables:
$$r_1=1 \qquad x_1=1 \qquad y_1=0.$$ 
Moreover, the equations $f_1=0$ and $\Omega = \w$ simplify to
$$u_1=0 \qquad  \fr12 v_1+\fr23(x_2v_2-y_2u_2) = \w$$
 eliminating two more variables, $u_1$ and $v_1$.  The remaining variables $r_{13},r_{23}$, $x_2,y_2,u_2,v_2$ are still constrained by the equations $f_2=f_3=f_4=0$, the energy equation and the equations for the squared mutual distances.  
  
At time $t_0$ we have $r_{13}^2 =(x_2+\fr12 )^2+ y_2^2$ and $r_{23}^2 =(x_2-\fr12 )^2+ y_2^2$ which allow us to eliminate $x_2$ and $y_2^2$.
\begin{equation}\label{eq_x2y2sq}
\begin{aligned}
x_2&= \fr12(r_{13}^2-r_{23}^2) \\
y_2^2&= \fr14(r_{13}+r_{23}+1)(r_{13}+r_{23}-1)(r_{13}-r_{23}+1)(r_{23}-r_{13}-1).
\end{aligned}
\end{equation}
Furthermore, if we introduce new velocity variables $w_{13}= \dot r_{13}$, $w_{23}=\dot r_{23}$ we have
$$r_{13}w_{13} = (x_2+\fr12)u_2+y_2v_2\qquad r_{23}w_{23} =(x_2-\fr12)u_2+y_2v_2.$$
These can be solved to find $u_2, v_2$ as functions of $r_{13},w_{13},r_{23},w_{23},y_2$.

The next step is to eliminate $y_2$ to obtain equations for the remaining, rotationally invariant,  variables.  Using the equation (\ref{eq_x2y2sq}) for $y_2^2$ we can reduce the numerators of the functions $f_2,f_3,f_4, H-h$ to functions involving only the first power of $y_2$.  Solving the simplest of these, $f_3=0$, for $y_2$ and substituting into the others gives three equations in $r_{13},w_{13},r_{23},w_{23}$ and substitution into the equation  (\ref{eq_x2y2sq}) for $y_2^2$  gives one more.  These equations, which are very complicated polynomials, will be denoted $g_i(r_{13},w_{13},r_{23},w_{23})=0$, $i=1,2,3,4$.  Their coefficients involve only integers and the parameters $h,\w$.  

Since the $g_i$ are function of the distances and their time derivatives, they are rotation-invariant.  We can use this rotation invariance to show that the equations $g_i=0$ apply not just at time $t=t_0$ but for every time.  Alternatively, we will show that these equations apply at a time $t=t_0$ even if $\g(t_0)$ does not satisfy the normalization requirement that $z_1(t_0)$ lies along the positive $x$-axis.  To see this suppose $\g(t_0)$ is an arbitrary point on a solution with $r_{12}=1$.  There is some fixed rotation of the plane, $R$,  such that $R\g(t_0)$ does lie along the axis.  Consider the curve $\tilde\g(t) = R\g(t)$.  By symmetry, this will represent another solution of the three-body problem with $r_{12}(t)=1$ which does satisfy the normalization condition at $t=t_0$.  It follows that the coordinates of $\tilde\g(t_0)$ must satisfy the equations $g_i=0$.  But since the equations $g_i$ are rotation-invariant, $0=g_i(\tilde\g(t_0))=g_i(R\g(t_0))=g_i(\g(t_0))$.

Because of the complexity of the equations, it's not possible to directly apply standard techniques of elimination theory such as resultants or Groebner bases.  Instead we will exploit a special structural feature of the equations to reduce the problem to two equations in the distances $r_{13}, r_{23}$ alone.  It happens that each of the equations $g_i$ is of degree two in the velocities $w_{13},w_{23}$ but none contain any terms of degree one.  Each can be written in the form
$$g_i = a_{i0} + a_{i1}w_{13}^2+a_{i2}w_{13}w_{23}+a_{i3}w_{23}^2 = 0$$
where the coefficient $a_{ij}$ are polynomials in $r_{13},r_{23}$.  In other words, they are linear combinations of the monomials $1, w_{13}^2, w_{13}w_{23}, w_{23}^2$.  It follows that the determinant of the $4\times 4$ matrix with entries $a_{ij}$ must vanish.  This gives a very large polynomial equation $G(r_{13},r_{23})=0$.  It follows that the distances $r_{13}(t),r_{23}(t)$ must satisfy $G(r_{13}(t),r_{23}(t))=0$ for all $t$.  To get another equation in two variables, first differentiate this equation to find
$$P(r_{13},r_{23})w_{13} + Q(r_{13},r_{23})w_{23} = 0$$ 
where $P,Q$ are the partial derivatives of $G$.  Then the equation $P^2w_{13}^2= Q^2w_{23}^2$ gives another equation in the same monomials.  The resulting $5\times 4$ matrix must have rank at most 4.  There will be five $4\times 4$ determinants $G_i(r_{13},r_{23})$, $i=1,\ldots,5$  obtained by deleting one of the rows. All of these determinants must vanish.

To finish the proof of the theorem, it remains to show that the five equation in two variables $G_i(r_{13},r_{23})=0$ have only finitely many solutions.  Again, the equations are too complicated for standard methods.  For example, the most complicated $G_i$ involves 4450 monomials including powers up to 105 in each variable.   In the next section, we will use Newton polygons to rule out a continuum of solutions, thereby completing the proof.

\section{Finiteness proof}
The proof involves the same ideas used in \cite{HamMoe} to prove finiteness of relative equilibria in the $4$-body problem and in \cite{Moe2} to prove Saari's conjecture.  Namely,  if there were infinitely many solutions to the equations $G_i(r_{13},r_{23})=0$ with $r_{ij}\in \C$,  there would be a solution in the form of nonconstant Puiseux series $r_{ij}(s)$ where $s$ is a complex parameter.  These series would be of the form 
$$r_{13}(s) = k_{13}s^a+\ldots\qquad r_{23}(s)=k_{23}s^b\ldots$$
where the lowest order terms have rational exponent vector $(a,b)\in \Q\times\Q$ and nonzero coefficients $k_{ij}\in \C$.  As explained in \cite{HamMoe} and elsewhere, the vector 
$(a,b)$ is an inner normal  of the Minkowski sum of the Newton polyhedra of the equations.  The coefficients $k_{13},k_{23}$ solve a reduced system of equations formed by taking, for each $G_i$, only those monomials $r_{13}^m r_{23}^n$ for which $am+bn$ is minimal.  For each $i$, this minimization determines a face of the Newton polytope.  If the face reduces to a vertex for some $i$ there will be no nonzero solutions for  $k_{13},k_{23}$.  

Using these ideas, we see that it suffices to check all of the faces of the Minkowski sum polytope to see if the corresponding reduced system of equations has a nonzero solution.  Since our equations involve only two variables, the Newton polytopes and their Minkowski sum are plane polygons.  The Newton polygon of $G_i$ is the convex hull of the set of all exponent vectors of the monomials appearing in $G_i$ with nonzero coefficients.  Using Mathematica, it's easy to find these exponent vectors and the convex hulls can be computed using the built-in command ConvexHullMesh[ ].  The Minkowski sum polygon can be found as the convex hull of the sum of the sets of exponents of the five $G_i$.  Although the coefficients of the $G_i$ depend on the parameters $h,\omega$, it turns out that the terms corresponding to the vertices of the Newton polygons are just nonzero integers, so the polygons are parameter-independent.

\begin{figure}[h]
\scalebox{0.6}{\includegraphics{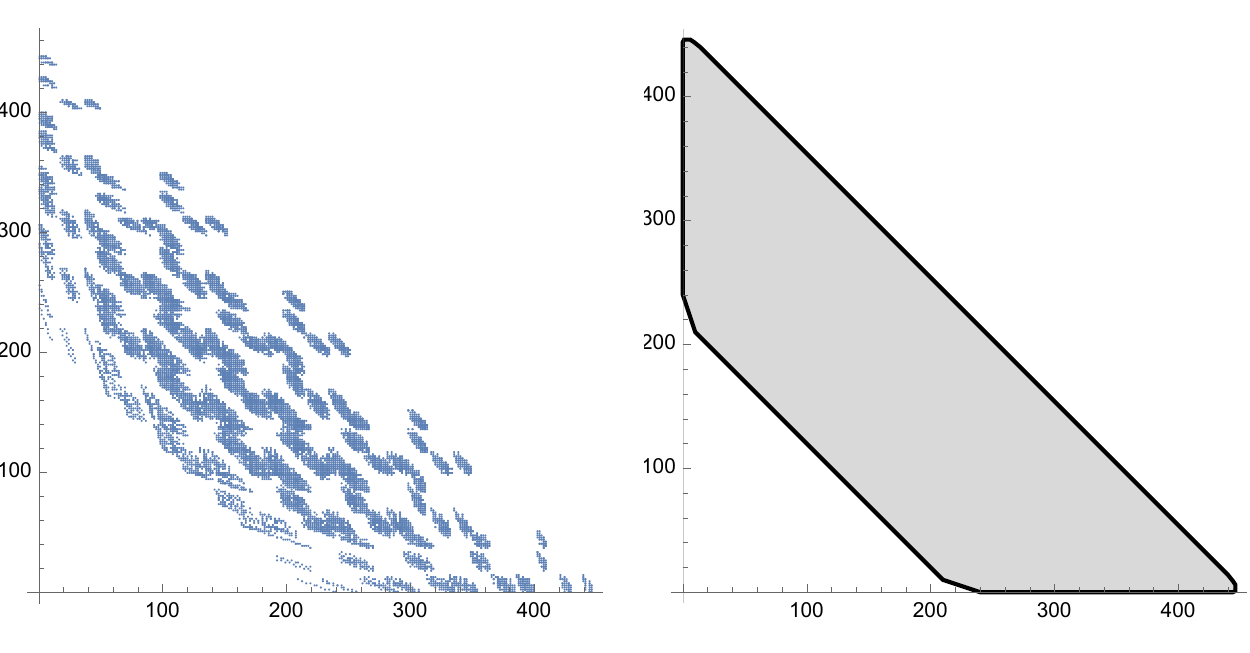}}
\caption{The Minkowski sum polygon of the Newton polygons of $G_1, \ldots, G_5$.  The sum of all the exponent vectors contains 16854 points (left) but the convex hull is quite simple (right).} \label{fig_MinkSum}
\end{figure}

Figure~\ref{fig_MinkSum} shows the sum of the five sets of exponents, which consists of 16854 integer points.  Fortunately, the convex hull is rather simple.  It's a polygon with 14 vertices and edges -- a tetradecagon.

As further explained in \cite{HamMoe} it suffices to check only some of the 14 edge, namely those for which the inner normal satisfies the inequality $a+b\ge0$, or similar.  These are the seven edges on the lower-left half of Figure~\ref{fig_MinkSum} which turn out to have inner normals
$$(2,-1)\quad (1,0)\quad (3,1)\quad (1,1)\quad (1,3)\quad (0,1)\quad (-1,2).$$
The reduced systems of some of them are quickly found to have no nonzero solutions.  For example, the reduced system corresponding to the inner normal $(3,1)$ consists of five simple equations, which, after dropping integer factors and factors of power of $r_{13}$ and $r_{23}$ contains the following:
$$3r_{13}+r_{23}^3=0\qquad 12636 r_{13}^3+9828 r_{13}r_{23}^3 + 1915 r_{23}^6=0\qquad \ldots$$
which have no nonzero common solutions.  Note that these equations are quasi-homogeneous in the sense that they are invariant under the scaling
$(k^3 r_{13},k^1 r_{23})$.  This permits us to add a normalization condition such as requiring $r_{13}r_{23}=1$.  Adding this equation to the reduced system and finding a Groebner basis is a way to automate the proof that there are no nonzero solution.

In this way, it is found that the only possibilities for the exponent vector of a Puiseux solution are $(1,1), (1,0), (0,1)$.  These need closer examination.  Consider the exponent vector $(1,1)$.   Both series must start with terms of the form $k_{ij}s$ with $k_{ij}\ne 0$.  In fact, by  series inversion, we can assume without loss of generality that one of the series is identically equal to $s$.  Setting $r_{23}=s$ and $r_{13}=k_{13}s$ where $k_{13}$ is a  Puiseux series with nonzero constant term, the lowest order terms of the five series
$G_i( k_{13}s,s)$ agree with the reduced system of vector $(1,1)$.  These five equations do have a nonzero common solution.  A Groebner basis consists of the single polynomial
$$P(k_{13})=k_{13}(1+k_{13})=0.$$
Thus the constant term is uniquely determined as $k_{13}=-1$.  Moreover we can use the implicit function theorem to argue that any Puiseux series solution
$k_{13}(s) = -s+\ldots$ must be a power series in $s$, that is, $k_{13}(s) = -s+u(s)$ where $u(s) = as^2+\ldots$  (rather than a series in fractional powers).  To see this, suppose that the lowest-order term of $G_i( k_{13}s,s)$ is of degree $d_i$ so that
 $$G_i( k_{13}s,s) = s^{d_i}F_i( k_{13}s,s)$$
 for some polynomial $F_i$ with $F_i(k_{13},0)\ne 0$.  The Groebner basis computations shows that there are polynomials $\phi_i(k_{13})$ such that 
 \begin{equation}\label{eq_F}
 F( k_{13}s,s)= \phi_1(k_{13})F_1( k_{13}s,s)+\ldots + \phi_5(k_{13})F_5( k_{13}s,s) = P(k_{13})+O(s).
 \end{equation}
 A nonzero Puiseux solution of $G_i( k_{13}s,s)=0$ must satisfy $F( k_{13}s,s)=0$.  Now $F(k_{13},0)=P(k_{13})=k_{13}(1+k_{13})$. Since $k_{13}=-1$ is a nondegenerate root,  the implicit function theorem guarantees that any solution $k_{13}(s)$ with $k_{13}(0)=-1$ is actually analytic, hence given by a series in integer powers.
 
 To show that no such series exists, we can now set $k_{13}(s) = -s+u(s)$ with $u(s)=a_2\, s^2 + a_3s^3+a_4s^4 \ldots$.  Substituting into each $G_i( k_{13}(s),s)$ and examining the  lowest order terms in $s$ shows that a solution is possible only if $a_2=0$ and furthermore the parameters $h, \w$ satisfy $h=-1+\fr32\w^2$.  Making these substitutions, the next lowest terms imply that $a_3=0$.  Making this substitution the next lowest terms show that there is no solution for $a_4$.  Namely, we arrive at an incompatible set of five equations equations including, for example,
 $$3a_4+36\w^2-8 = 0\qquad 189a_4+3438\w^2-704=0.$$
 This finally shows that no nonzero Puiseux solution with exponent vector $(1,1)$ exists.
 
\begin{figure}[h]
\scalebox{0.4}{\includegraphics{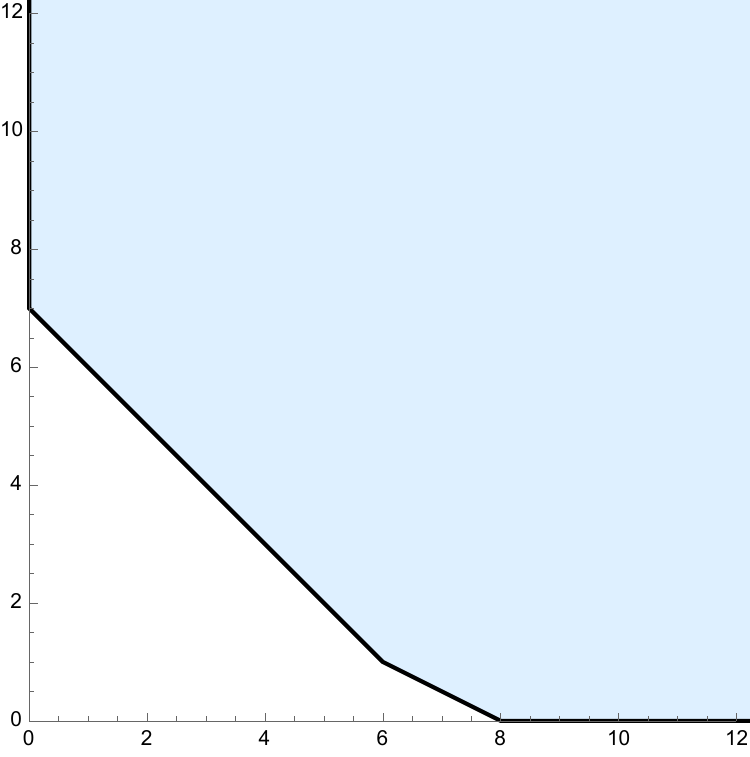}}
\caption{The Newton polygons of $F_5$, the simplest of the five equations $F_i(s,u)$. The slopes of the line segments are $-1, -\fr12$} \label{fig_Newton}
\end{figure}

 We can try a similar approach for the exponent vectors $(0,1)$.  Setting $r_{23}=s$, the Puiseux series for $r_{13}=k_{13}(s)$ would begin at order $0$.  
 A Groebner basis for the five lowest order terms is
 $$P(k_{13}) = k_{13}^{21}(k_{13}-1)^7(k_{13}+1)^7.$$
 It follows that the constant term must be $k_{13}= \pm 1$.  But this time, the roots of $P$ are degenerate so we cannot immediately assume a power series.  Consider the case
 $k_{13}(s) = 1 +u(s)$ where $u(s)$ is a Puiseux series with $u(0)=0$.  Substitution gives five equations $F_i(s,u)=G_i(1+u,s)$ with $F_i(0,0)=0$ and we seek a solution $u(s)$ with $u(0)=0$.  This is the classic case for the Newton diagram method for finding branches of an algebraic curve.  The lower convex hull of the set of exponents of each equation is bounded by a polygonal curve of  line segments with negative slopes $m$ connecting the coordinate axes.  The  lowest-order term of any Puiseux solution $u(s)$ must be $d=-\fr{1}{m}$ for some slope $m$ (since $(1,-\fr1m)$ is an inner normal).  Figure~\ref{fig_Newton} shows the Newton diagram for the simplest of the five equations which consists of two line segments of slopes $-1,-2$.  So the lowest-order term in any Puiseux solution must be of the form $as$ or $as^2$.

Substituting $r_{13}(s) = 1+as$ where $a$ is a Puiseux series and finding the lowest-order term shows that there is no consistent solution for the constant term of $a$.  Similary, trying 
$r_{13}(s) = 1+as^2$ also leads to inconsistency.  Thus there is no Puiseux solution of $G_i=0$ with exponent vector $(0,1)$ and constant term $k_{13}=1$.  The argument can be repeated for the case $k_{13}=-1$.   Ruling out the last exponent vector $(1,0)$ follows the same pattern.   Finally we see that the requirement that $r_{12}(t)=1$ constrains
$r_{13},r_{23}$ to a finite set which completes the proof of  Theorem~\ref{th_nopartialrigidity}.

\end{document}